# Noisy data clusters are hollow

François Léonard

Hydro-Québec, IREQ, 1800 boul. Lionel-Boulet, Varennes (Québec), Canada J3X 1S1

**Abstract**

A new vision in multidimensional statistics is proposed impacting several areas of application. In these applications, a set of noisy measurements characterizing the repeatable response of a process is known as a realization and can be seen as a single point in $\mathbb{R}^N$. The projections of this point on the $N$ axes correspond to the $N$ measurements. The contemporary vision of a diffuse cloud of realizations distributed in $\mathbb{R}^N$ is replaced by a cloud in the shape of a shell surrounding a topological manifold. This manifold corresponds to the process's stabilized-response domain observed without the measurement noise. The measurement noise, which accumulates over several dimensions, distances each realization from the manifold. The probability density function (PDF) of the realization-to-manifold distance creates the shell. Considering the central limit theorem as the number of dimensions increases, the PDF tends toward the normal distribution $N(\mu, \sigma^2)$ where $\mu$ fixes the center shell location and $\sigma$ fixes the shell thickness. In the proposed vision, the likelihood of a realization is a function of the realization-to-shell distance rather than the realization-to-manifold distance. The demonstration begins with the work of Claude Shannon followed by the introduction of the shell manifold and ends with practical applications to monitoring equipment.

## 1. Introduction

This paper focuses on a link between the measurement noise of a repeatable process and hollow clusters in $\mathbb{R}^N$ and proposes a method to exploit this link in different applications. The proposed link has been only partially revealed by papers about fuzzy clustering, shell clustering, noise clusters, hollow clusters in imaging, vacuum clusters, data topology, topological noise and persistent homology.

In a generalization of the fuzzy k-means clustering algorithm, Coray [1] suggested the use of non-linear circular prototypes for clusters. Some shell-clustering algorithms and their application to boundary detection and surface approximation have been proposed by Krishnapuram et al. [2] with the concept of expected thickness of the shell clusters. Some authors [3,4] propose clustering tools adapted to noisy data. For example, a noise cluster prototype is introduced in the hope that all the noisy points can be part of the noise cluster [5,6]. Hollow clusters have been observed in imaging [7]. In pattern-recognition techniques, vacuum shell detection is used in clustering methods [8].

The topic of data topology is linked to many theoretical aspects and application fields. Hollow shapes can be detected by homology, which is a computable topological invariant of shapes. Niyogi et al. [9] show that if the data is noisy, in the sense that it is drawn from a probability distribution concentrated around a manifold, the homology of that manifold can still be computed with high confidence from that data. In topology, the noisy data is handled more systematically by the persistent homology methods developed recently by several authors [10-12]. Pekelis and Holmes [13] present statistical approaches to calculating the homology of data and testing whether data lie on a closed manifold with a hollow interior. The reader is referred to Carlsson [14] for a comprehensive review of this topic.

The link between measurement noise and hollow clusters starts with the definition of repeatable process. A repeatable process offers an opportunity for pattern recognition, monitoring, trending, diagnosis and prognosis. The time dimension is implicit in the context of monitoring but not necessary in pattern recognition. Repeatable processes are found in many fields: medical, biological, chemical, mechanical, electrical, quantum mechanics, social, stock market, manufacturing, communication and imaging. A process can be defined as a transfer function having internal states, external inputs (independent variables) and external measured responses (dependent variables).

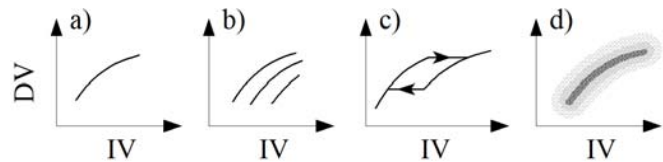

**Figure 1:** 2D illustration of the response function, a dependent variable (DV) function of an independent variable (IV), for a) a single internal state, b) several internal states, c) a hysteresis internal state and d) an additive effect of numerous internal states.

A process can be observed from the measurements of independent and dependent variables. Measurements contain both information and noise. With regard to *information,* in some specific cases, e.g. communications, the information set is known and defined by a finite set of discrete values. In more common cases, the information is unknown and can be defined by a multidimensional function corresponding to a unique internal state (Fig. 1a). As for *noise*, this is a zero-mean random process, correlated or not between measurements. Noise can be seen as the sum of numerous sources of irrecoverable information. Hidden state processes with a small impact contribute to noise. When the time dimension is implicit, a *realization* is a set of close time measurements recorded on a same channel or on many channels. In practical terms, the operating condition of the process is time-variant and, as in photography, the realization is impaired by a blur when the time spread of the measurement set is too large. In the following development, a realization corresponds to a picture of a stabilized process and can be taken discontinuously with an irregular sampling rate. Note that the "time" dimension can be replaced by the appropriate dimension related to the context (e.g. spatial). The measurement set corresponding to a realization can share a common engineering unit or not. The normalization of measurements having different engineering units is an item of great concern: the selected normalization method contributes to the success or failure of the application algorithm (monitoring, pattern recognition…). Turning lastly to *internal*

*states,* a response function can contain more than one internal state which can be visible or not from the dependent variables and can be determined by a sequence of independent variable values (Markov chain) or hidden. Figure 1 illustrates the response function corresponding to a single, several and multiple internal states. A hidden state can be defined as a state not linked to the input variables with causality. In the latter case, an unmeasured independent variable which contributes inputs to the process can be transposed as a hidden state of the process. The additive effect of numerous small contributions of unmeasured independent variables and hidden states to the process yields a diffused response function (Fig. 1d).

In multidimensional statistics, a realization corresponds to a single point in $\mathbb{R}^N$. A set of realizations corresponding to the same state process yields a cluster in $\mathbb{R}^N$. This paper starts with the Shannon vision of the cluster constellation corresponding to the received signals: a stack of *n*-spheres with *n=N-1*, each *n*-sphere centered on its corresponding transmitted symbol before the addition of noise. From this vision of multiple *n*-sphere clusters associated with already known symbols having no dimension (a point in $\mathbb{R}^N$ for each state process), we will jump to a single cluster generated by an unknown multivariate one-state process. In both cases, the clusters are hollow.

## 2. Back to Shannon

Hollow clusters of measured data were first introduced by Claude Shannon [15] in the mid 20th century. His demonstration of the limit capacity of a communication channel is based on the limit of the stacking of *n*-spheres in a $\mathbb{R}^N$ volume delimited by a maximum amplitude signal. The *n*-sphere radius corresponds to the average length of the noise vector added to a transmitted symbol. The error rate appears as a function of the overlap between the shells, each shell corresponding to a transmitted symbol.

In Shannon's demonstration, a communication channel comprises an information source (the already known symbol set), a transmitter, a noise source, a receiver and a destination. Consider only one symbol

$$\mathbf{x} \equiv [x_1, \ldots, x_N] \tag{1}$$

of the set transmitted using *N* dimensions. The corresponding received signal is

$$\mathbf{y} = [x_1 + \varepsilon_1, \ldots, x_1 + \varepsilon_n, \ldots, x_N + \varepsilon_N] = \mathbf{x} + \boldsymbol{\varepsilon} \tag{2}$$

where $\varepsilon_n$ is the noise added for the dimension "*n*" and $\boldsymbol{\varepsilon}$ the noise vector added to the transmitted symbol in $\mathbb{R}^N$. The average distance from the signal received to the symbol transmitted

$$\mu = \mathrm{E}\big[\|\mathbf{y} - \mathbf{x}\|\big] \tag{3}$$

corresponds to the radius of the *n*-sphere and becomes

$$\mu = \mathrm{E}\left[\sqrt{\sum_{n=1}^{N} \varepsilon_n^2}\right] \tag{4}$$

considering the use of the Euclidian metric. The measurement dispersion

$$\sigma^2 = \mathrm{E}\left[\left(\|\mathbf{y} - \mathbf{x}\| - \mu\right)^2\right] \tag{5}$$

is called the hardness of the $n$-sphere. The increasing $\mu/\sigma$ ratio with the increasing number of dimensions is called the "sphere-hardening phenomenon" in the communication theory [16]. In the present mathematical development, $\sigma$ is defined as the half-thickness of the shell. Indeed, when the added noise $\{\varepsilon_n\}$ is a set of $N$ zero-mean Gaussian PDF independent random variables, the PDF of a realization is the normal distribution $\mathrm{N}(\mu, \sigma^2)$. However, when the added noise $\{\varepsilon_n\}$ is a set of $N$ zero-mean random variables uncorrelated, the central limit theorem stipulates that the PDF tends also toward the normal distribution $\mathrm{N}(\mu, \sigma^2)$ with the increasing number of dimensions.

In a practical application where $\mathbf{x}$ is unknown and $\mathbf{y}_m$ is the realization "$m$", eqs. 4 and 5 are replaced by the estimates

$$\hat{\mu}_M = \frac{1}{M} \sum_{m=1}^{M} \|\mathbf{y}_m - \overline{\mathbf{y}}_M\| \tag{6}$$

and

$$\hat{\sigma}_M^2 = \frac{1}{M} \sum_{m=1}^{M} \left(\|\mathbf{y}_m - \overline{\mathbf{y}}_M\| - \hat{\mu}_M\right)^2 \tag{7}$$

where

$$\overline{\mathbf{y}}_M = \frac{1}{M} \sum_{m=1}^{M} \mathbf{y}_m \tag{8}$$

is the estimated value of $\mathbf{x}$ calculated from $M$ realizations.

The PDF of a realization is defined by

$$\mathrm{N}\left(\hat{\mu}_M, \hat{\sigma}_M^2\right) \text{ when } \{\mathbf{y}_m : m \in [1, M]\} \text{ and by} \tag{9}$$

$$\mathrm{N}\left(\hat{\mu}_M + \varepsilon_{\mu_M}, \hat{\sigma}_M^2\right) \text{ when } \{\mathbf{y}_m : m \notin [1, M]\} \text{ where} \tag{10}$$

$$\varepsilon_{\mu_M} = \hat{\mu}_M \cdot \sqrt{\frac{1}{M(M-1)}} \tag{11}$$

is an additional distance explained by a new realization not taken into account in the $\overline{\mathbf{y}}_M$ estimation. The case of a realization as an average number of realizations and that of the comparison of two estimations $\overline{\mathbf{y}}_{M_1}$ and $\overline{\mathbf{y}}_{M_2}$ calculated for different sets of realizations are discussed in [17].

## 3. Shell manifold

Let

$$x_n = f_n(w_1, w_2, \ldots, w_l, \ldots, w_L), \ n \in [1, N], \qquad (12)$$

be the *N*-coordinates system of the manifold *X* in $\mathbb{R}^N$ defined by *L* variables $w_l$. Note that *L*=0 for the specific case of the *n*-sphere. In the general context of an application, the *N* dependent variables are the response of the process and the *L* independent variables are the process inputs.

Consider the $m^{\text{th}}$ realization

$$\mathbf{y}_m = \left[ x_{1,m} + \varepsilon_{1,m}, \ \ldots, x_{n,m} + \varepsilon_{n,m}, \ \ldots, x_{N,m} + \varepsilon_{N,m} \right] = \mathbf{x}_m + \boldsymbol{\varepsilon}_m \qquad (13)$$

of the realization set *Y* with

$$\mathbf{y}_m = \left[ y_{1,m}, \ldots, y_{n,m}, \ldots, y_{N,m} \right], \ \mathbf{y}_m \in Y, \qquad (14)$$

and where $\varepsilon_{n,m}$ is the noise added for the dimension *n* and $\boldsymbol{\varepsilon}_m$ is the noise vector in $\mathbb{R}^N$ which contributes to the realization *m*.

Consider the cluster corresponding to repeatable realizations $\mathbf{y}_m = \mathbf{x}_{\text{point}} + \boldsymbol{\varepsilon}_m$ for a unique point $\mathbf{x}_{\text{point}}$ of the manifold, the average distance from the realization to the point is

$$\mu_{\text{point}} = \mathrm{E}\left[ \|\mathbf{y}_m - \mathbf{x}_{\text{point}}\| \right] \qquad (15)$$

and corresponds to the radius of the *n*-sphere eq. 3.

In the case of realizations set *Y* corresponding to different points on the manifold, considering the special case of a equal amplitude additive Gaussian noise $\mathrm{N}(0, \varepsilon_0^2)$ for all dimensions and a local small surface patch S having a negligible curvature with respect to the noise vector length with *N-L*>>1, we can apply some axis rotations to $\mathbb{R}^N$ in order to align the local *L*-surface patch S of the manifold in the bearing of *L* axis of the new *N* coordinates system of $\mathbb{R}^N$. The mean distance of the realization normal to the surface

$$\mu_{\perp \mathrm{S}}\big|_{\varepsilon_n = \varepsilon_0} = \mathrm{E}\left[ \sqrt{\sum_{n=1}^{N-L} \varepsilon_n^2} \right]_{\varepsilon_n = \varepsilon_0} = \varepsilon_0 \sqrt{N-L} = \mu_{\text{point}} \sqrt{\frac{N-L}{N}} \qquad (16)$$

is shorter than the noise vector since only *N-L* coordinates contribute to this distance. Conversely the corresponding variance

$$\sigma_{\perp \mathrm{S}}^2\big|_{\varepsilon_n = \varepsilon_0} = \mathrm{E}\left[ \left( \sqrt{\sum_{n=1}^{N-L} \varepsilon_n^2} - \mu_{\perp \mathrm{S}} \right)^2 \right]_{\varepsilon_n = \varepsilon_0} = \varepsilon_0^2 \frac{N}{2(N-L)} = \mu_{\text{point}}^2 \frac{1}{2(N-L)}. \qquad (17)$$

increases when the number of manifold dimensions *L* increases. This simplest case of a Gaussian noise having a uniform amplitude $\varepsilon_0$ on the *N* dimensions and independent of the manifold coordinates yields a hollow cluster of realizations where the PDF of a realization is $N\left(\varepsilon_0\sqrt{N-L},\varepsilon_0^2\cdot N/2(N-L)\right)$. In the general case of the realizations set **Y** built from dimensions having different noise sources, we conclude that locally $\mu_{\perp S} < \mu_{point}$ in the vicinity of a surface patch S. For the general case of *M* realizations locally available near the surface patch S, we define the average realization-to-manifold distance as

$$\mu_{\perp S} = \frac{1}{M}\sum_{m=1}^{M}\|\mathbf{y}_m - \hat{\mathbf{x}}_m\| \tag{18}$$

where $\hat{\mathbf{x}}_m$ is the nearest surface point of the estimated manifold in the vicinity of the surface patch S with respect to the realization $\mathbf{y}_m$. This average distance fixes the shell position with respect to the manifold. Since *N>L*, the shell surrounds the manifold. The estimation of realization dispersion relative to the surface

$$\sigma_{\perp S}^2 = \frac{1}{M}\sum_{m=1}^{M}\left(\|\mathbf{y}_m - \hat{\mathbf{x}}_m\| - \mu_{\perp S}\right)^2 \tag{19}$$

is called the hardness of the shell. The $\sigma_{\perp S}$ value is defined as the half-shell thickness.

Note that the average distance between two realizations $\mathbf{y}_i$ and $\mathbf{y}_j$ corresponding to the same set of independent variables (same location on the manifold) is

$$E\left[\|\mathbf{y}_i - \mathbf{y}_j\|\right] = \mu_{\perp S}\cdot\sqrt{\frac{2N}{N-L}} = \varepsilon_0\cdot\sqrt{2N} \tag{20}$$

for the simplest case of a Gaussian noise having a uniform amplitude $\varepsilon_0$ on the *N* dimensions and independent of the manifold coordinates. The nearest neighbor distance of a realization is on average less than the previous average distance (eq. 20) and diminishes with the increasing measurement population. When the noise is a function of the manifold coordinates, the shell position and hardness also vary with the manifold coordinates. In the more general case of a diffuse manifold (Fig. 1d), which is not defined by eq. 12, the PDF of the realizations corresponds to the multidimensional convolution of the manifold density with the $N\left(\mu_{\perp S},\sigma_{\perp S}^2\right)$ function. The cluster appears hollow only if the total noise amplitude is significantly greater than the manifold width. In the case of multiple internal states, the resulting PDFs of the realization can be superimposed. To sum up, Fig. 2 shows an artistic 3D illustration of a hollow cluster and the PDF of realizations when the manifold is known (b) and when it is estimated (c). When the manifold is estimated, the quadratic form $\left(\sigma_\mathbf{M}^2 + \sigma_{\perp S}^2\right)$ of the combined dispersion is used to estimate the likelihood of a realization where $\sigma_\mathbf{M}$ is the manifold dispersion estimated in the vicinity of the realization.

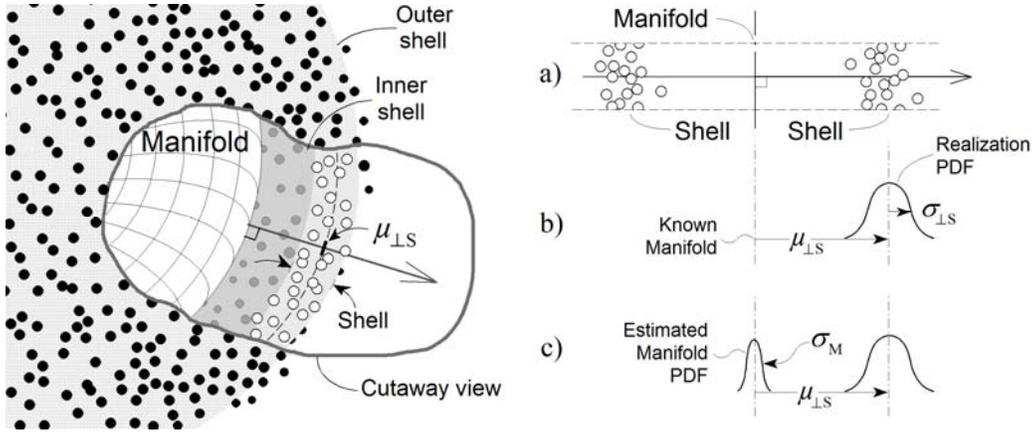

**Figure 2:** Cutaway view of a cluster showing a manifold (left). Realizations on the outside shell are in black, on the inside in grey and in the cut view in white. On the right of the figure: a) a cutaway view through the manifold surface, b) the corresponding PDF for a known manifold and c) the general case of an unknown manifold estimated from realizations averaging.

## 4. A sample application with a generic comparison algorithm

The algorithm compares two vector quantities X and Y in $\mathbb{R}^N$. The vector may be a realization or an estimate (bin average or manifold interpolation). The average of the distance ||X-Y|| is assumed constant until the repeatable response of the process is no longer there. In the following algorithm sequence, the "match" threshold level is set at four times the standard deviation $shellvar^{1/2}$:

1. $m = 0$; 'Number of processed comparisons'
2. $shelldist = 0$ : 'Mean manifold-to-shell distance'
3. $shellvar = 0$ ; 'Variance of calculated distances to the mean distance'
4. $match = .true.$ ;
5. Wait, calculate or find in memory the X and Y vectors;
6. $d = ||X-Y||$ ; 'Calculate the distance'
7. If $m=0$ Then Goto 10 ; Endif ; If $m=1$ Then Goto 9 ; Endif;
8. If $|d - shelldist| > 4 \cdot shellvar^{1/2}$ then $match = .false.$ ; else $match = .true.$ ; Endif ;
9. $shellvar = ((d - shelldist)^2 + (m-1) \cdot shellvar) / m$ ;
10. $shelldist = (d + m \cdot shelldist) / (m+1)$ ;
11. $m = m+1$ ;
12. goto 5 ;

Lines 8 to 10 are the kernel of the application algorithm and contain the essential part of the innovation. Lines 9 and 10 can be replaced by moving averages as:

9. $shellvar = \left(alpha \cdot (d - shelldist)^2 + (1-alpha) \cdot shellvar\right) / alphaweight$
10. $shelldist = \left(alpha \cdot d + (1-alpha) \cdot shelldist\right) / alphaweight$

where $alpha \in ]0,1[$ is the smoothing factor of the moving average with a weighting factor $alphaweight = alpha + (1 - alpha) \cdot alphaweight$ calculated before line 9 and with $alphaweight = 0$ at the initialization step. The contributing variances *Ex* and *Ey* of the estimates X and Y can be included in the line 8 as:

8. If $|d - shelldist| > 4 \cdot (shellvar + Ex^2 + Ey^2)^{1/2}$ then *match*=. false. ; else *match*=. true.

These errors correspond to the estimation error vector length in $\mathbb{R}^N$. If X or Y is a realization, no error is estimated for it. For the case of a clustering followed by an interpolation (Fig. 3), the *Ex* or *Ey* errors correspond to the manifold dispersion $\sigma_\mathbf{M}$. This dispersion can be estimated from error propagation in the clustering process taking account of the weight of each cluster in the interpolation and the distance between an interpolated point on the manifold with the surrounding clusters. More computationally expensive, the manifold dispersion estimation can be obtained by bootstrapping [18]. Whatever the method used, the computational complexity deployed in these error estimations substantially exceeds that needed for process lines 6-11.

Lines 9 and 10 allow the detection of an abrupt distance change between X and Y. To detect a slow trend, the *shelldist* and *shellvar* values must be frozen when the population *m* is large enough to be representative of the average initial response. As a coding option, the *shelldist* and *shellvar* values can be estimated considering the new realization only when *match* is true.

## 5. Practical monitoring applications

A monitoring application software includes measuring, database management, numerical signal processing, estimation, comparison and alarm management. Signal processing transforms the raw measurements into a distilled form in a representation space (envelope, correlogram, cyclogram, Fourier, time-frequency, scale-frequency, high-order spectrum…) where the signal-to-noise ratio appears at its maximum in the context of the application task. The estimation builds the initial history of a process and refreshes the actual process average response. With the "match" logic case replaced by a "no-alarm" logic case (lines 4 and 8), the previous generic comparison algorithm can be used for monitoring. Usually two comparisons are done, one between the initial history and the actualized average response for slow trend detection, and the second between the actualized average response and the last realization for the detection of a fast growing defect [19].

In vibro-acoustic monitoring of commutating equipment, the signal envelope is compared to a reference envelope built from data history [19]. The equipment response is found in the *N* time samples of the vibro-acoustic envelope [20], typically 1000 samples. The *n*-sphere method described in eqs. 1-11 has been applied with limited success since the equipment response is significantly a function of the temperature. The binning of the reference envelope as a function of a set of temperature ranges was a first step toward increasing the accuracy of the method. The trajectory manifold (*L*=1) is modeled by a set of bins, each bin containing an average location in $\mathbb{R}^N$, an average temperature and the corresponding population. In a second step, we added a krigging algorithm [21-23] to interpolate the reference value from the bins. The last realization is then compared to the reference interpolated at the same temperature as the realization. The bin values, populations, temperatures and dispersions are used as krigging inputs. Each of the 1000 time samples is interpolated.

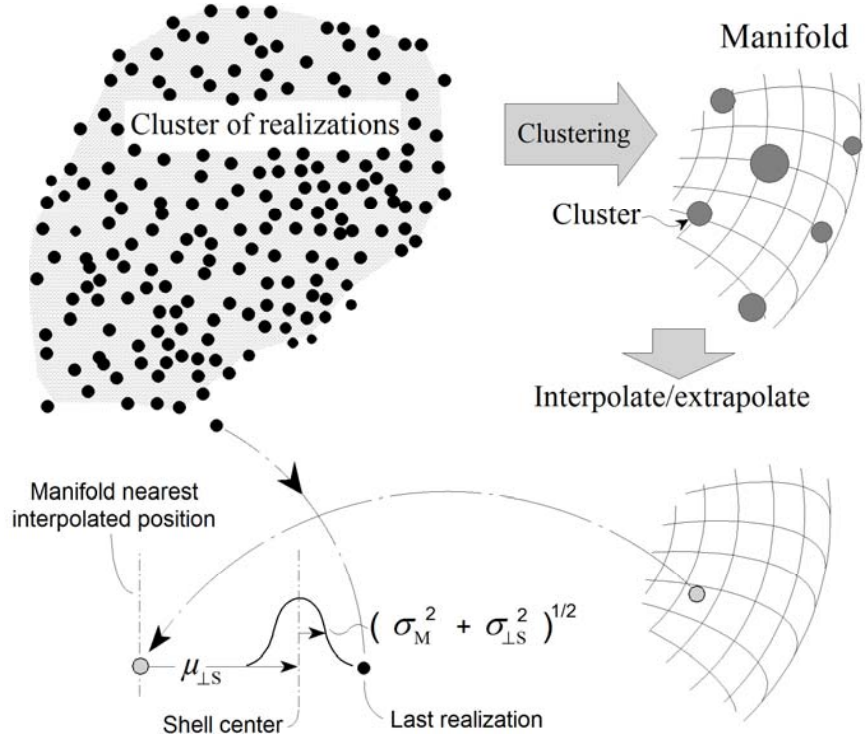

**Figure 3:** Sequence of operations in estimating the likelihood of the last realization. Realizations are clustered in order to estimate the inner manifold. An interpolation of this manifold is then compared to the last realization. In the upper right illustration, the cluster disk is proportional in size to the cluster population. Clustering estimates the values $\mu_{\perp S}$ and $\sigma_{\perp S}$, while krigging estimates the interpolation and its uncertainty $\sigma_M$.

In monitoring a hydroelectric generating unit [24,25], the operating condition is determined by many independent variables: gate opening, water levels, rotor excitation current, cooling water and air temperatures,… If we include vibrational data (peak-peak, RMS, mean and spectral line values) as part of the dependent variables, the total number of measured variables reaches over 100. In the proposed method, a dynamic clustering estimates the manifold corresponding to the behavior of the unit. The clustering in $\mathbb{R}^N$ is a dynamic clustering inspired from the one described in [26,27] and described in the addendum. Since the measured variables are expressed in different engineering units, we added a dynamic normalization of variables as a function of the maximum observed range of each variable. The clustering is usually applied using the information from the $L$ independent variables (see *mask* = 1 in the addendum). The clustering can also use the information about some dependent variables in order to distinguish between a limited number of internal states (Fig. 1b) but this feature has not been tested. The average value of each dependent variable and the average square value are calculated for each cluster.

When enough clusters with enough populations are found from the data history, krigging is used to interpolate the $N$ dependent variables averaged in clusters. Krigging uses the cluster values, dispersions and populations to estimate the interpolated (or extrapolated) value and its accuracy. The interpolated value and its accuracy $\sigma_M$ are both used in the comparison with the last realization.

For both applications, the accuracy of the interpolation is added quadratically (Fig. 3, lower) to the shell realization dispersion to fix the alarm bound (typically ±4 $\sigma$). Before the comparison with the alarm bound, the mean distance $\mu_{\perp S}$ is subtracted from the distance between the interpolation and the last realization. Note that every change in the equipment response pushes the resulting distance to increase toward the positive bound (+4 $\sigma$). The newly developed method dramatically reduces the rate of false alarms (type 1 error) for a same missing-defect rate (type 2 error).

Clustering (or binning) is an unavoidable part of the method. "Clustering should be thought of as the statistical counterpart to the geometric construction of the path-connected components of a space, which is the fundamental building block upon which algebraic topology is based" [1]. In the proposed applications, the clustering is done function of the $L$ independent variables assuming a single continuous and monotonous manifold. Moreover, the measurement noise of independent variables is not considered: this noise is implicitly transferred to dependent variables.

## 6. Conclusion

For a repeatable response process, data clusters built from noisy data are hollow when the noise amplitude is greater than the dispersion generated by the hidden internal states of the process. When the clusters are built from repeatable realizations, clustering and interpolation can be used to estimate the manifold corresponding to denoised information. The probability density function of realizations generates a shell enveloping the manifold and the corresponding model can be estimated from realizations. Big data subsets related to measurements are shell manifolds.

The proposed vision suggests a new way to process information about noise in multidimensional statistics. The PDF contribution of a single measurement channel in the context of a multichannel measurement has been explained and exploited. The likelihood of a realization appears as function of the realization-to-shell distance rather than the realization-to-manifold distance.

Many of the current tools of multivariate analysis are affected by this vison since they do not target hollow clusters or shell clusters. For each of these tools, the question is: can we find an algorithm improvement to adapt the tool to noisy data clusters? This new proposed vision is a call to everyone looking for new tools adapted to noisy data clusters including large multidimensional datasets.

# Addendum – Clustering algorithm

The normalization of input variables having different engineering units is not addressed in the following description of a clustering algorithm. Dynamic clustering is a compromise between optimality and computational complexity. In the proposed clustering algorithm, the final number of clusters *kmax* is predetermined and fixes the memory dimension. The cluster centroid locations [28] are adjusted dynamically at each step considering the last clustering result and the new realization. The decision rule is based on a Euclidian distance calculated from a subset of available dimensions. In a realization, *L* measurements correspond to independent variables and *N* others correspond to dependent variables. The function

$$mask(i) = \begin{cases} 1 : \text{independent variable} \\ 0 : \text{dependent variable} \end{cases} \quad \text{with } i \in [1, N+L] \tag{A-1}$$

defines which subset is used in the distance calculation

$$\left\| \mathbf{C}(k) - \mathbf{X} \right\|_{masked} = \left( \sum_{i=1}^{N+L} mask(i) \cdot \left( C_i(k) - X_i \right)^2 \right)^{1/2}. \tag{A-2}$$

The average realization-to-cluster centroid variance *cvar(k)* is calculated for each cluster based on a masked distance calculation. The square root of this variance is averaged in order to share a unique shell distance value, *shelldist*, for all clusters. The *shelldist* value adjusts the cluster-to-cluster sampling distance on the manifold dynamically. Note that the increasing *kmax* value reduces the sampling distance and, also, the average cluster population as well as the accuracy of the centroid estimate. The number of clusters *kmax* fixes the tradeoff between the sampling distance and the cluster population. Clustering is successful when the clusters are uniformly distributed on the manifold. For the same realization set, different sampling orders of realizations yield similar results.

In order to reduce the computational expense, the maximum allowed distance *dmax* is defined as $cdist \cdot shelldist$ and corresponds to the decision level to merge a realization to the nearest cluster without the computational task of considering a possible fusion of two clusters. One cluster fusion $k \to k-1$ decision implies the use of a *k*-Nearest Neighbors (*k*-NN) algorithm [29,30] realizing $k \cdot (k-1)/2$ cluster-to-cluster comparisons without optimization. *cdist* is set typically at about 1 to 3, a function of *L* and the desired accuracy. Sharing the same *dmax* value for all clusters forces the cluster variances to converge into a similar range which corresponds to a uniform sampling of the explored part of the manifold. The algorithm stores the first *kmax* realizations as new clusters and starts the subsequent clustering. Considering X as a new realization, **C**(*k*) as a matrix of the *kmax* cluster centroids, **V**(*k*) as a matrix of cluster dispersion, *p(k)* as a vector of cluster populations and *shelldist* as the average shell cluster distance, the algorithm steps are:

1. *kcount* =0; *dmax* =0; 'Realization counter and *dmax* initialization'
2. *p*(*k*:1 to *kmax*)) =0; *cvar* (*k*:1 to *kmax*))=0; 'Population and variance initialization'
3. **C**(*k*:1 to *kmax*)) =0; 'Cluster centroid initialization'
4. **V**(*k*:1 to *kmax*)) =0; 'Cluster variance initialization'
5. Read in the memory or wait for the new realization X; *kcount* = *kcount* +1 ;
6. If *kcount* < *kmax*  Then $\mathbf{C}(k) = \mathbf{X}$ ; goto 5 'Initialization in progress'; Endif ;

7. $sdrc = \min_k \{\|\mathbf{C}(k) - \mathbf{X}\|_{masked}\}$ ; 'Finds the Smallest Distance Realization-to-Cluster and the corresponding cluster $k$'
8. If $sdrc > dmax$ Then
    'Calculates the cluster-to-cluster distances triangle matrix and finds $dscc$, the smallest distance, and the corresponding cluster couple $(j,k)$';
      If $sdrc > sdcc$ Then 'merge the two corresponding clusters'
      $\mathbf{C}_{fusion} = (p(j) \cdot \mathbf{C}(j) + p(k) \cdot \mathbf{C}(k))/(p(j) + p(k))$;
      $\mathbf{V}(k) = (p(j) \cdot \mathbf{V}(j) + p(k) \cdot \mathbf{V}(k))/(p(j) + p(k))$;
      $cvar(j) = p(j) \cdot (cvar(j) + \|\mathbf{C}(j) - \mathbf{C}_{fusion}\|^2_{masked})$
      $\quad\quad\quad + p(k) \cdot (cvar(k) + \|\mathbf{C}(k) - \mathbf{C}_{fusion}\|^2_{masked})$;
      $\mathbf{C}(k) = \mathbf{C}_{fusion}$; $\mathbf{C}(j) = \mathbf{X}$; goto 10;
      Endif;
    Endif;
9. 'Merges the realization to the nearest cluster'
    $\mathbf{C}(k) = (\mathbf{X} + p(k) \cdot \mathbf{C}(k))/(p(k) + 1)$;
    $\mathbf{V}(k) = ((\mathbf{C}(k) - \mathbf{X}) \circ (\mathbf{C}(k) - \mathbf{X}) + p(k) \cdot \mathbf{V}(k))/(p(k) + 1)$; '$\circ$: Hadamard product'
    $cvar(k) = (p(k) \cdot cvar(k) + sdrc^2)/(p(k) + 1)$; $p(k) = p(k) + 1$;
10. $shelldist = (\sum_k p(k) \cdot cvar(k) / \sum_k p(k))^{1/2}$ ; '$dmax$ calculation'
    $dmax = cdist \cdot shelldist$ ; goto 5;

At the beginning of the clustering process, the *shelldist* is underestimated (line 9) since every new cluster has zero variance. An underestimated value of *dmax* increases the probability of cluster inter-distance matrix calculation, increasing the accuracy and also the computational expense. When a second member $\mathbf{X}_j$ is merged with a first member $\mathbf{X}_k$ of a new cluster, the variance $cvar(k) = 0.5 \cdot \|\mathbf{X}_i - \mathbf{X}_k\|^2_{masked}$ has its expectation lower bounded $\mathrm{E}[cvar(k)] \geq \mu$ considering eq. 20. As demonstrated in monitoring a hydroelectric generating unit in Fig. 4, the algorithm converges to minimize calls to cluster-to-cluster matrix calculation. Note the "0" count during initialization and the asymptotic pattern of the curve which converges to add new calls only when the sampling touches an unexplored part of the manifold.

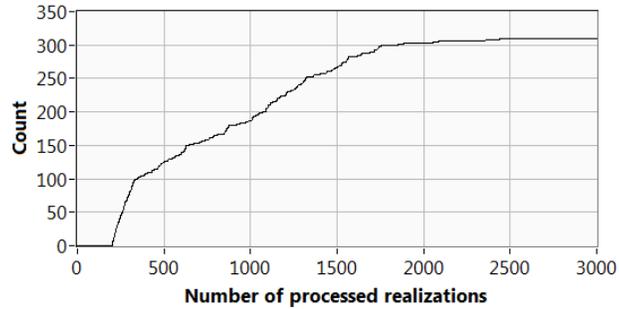

**Figure 4:** Counts of cluster-to-cluster matrix calculation as a function of the number of processed realizations. In this example, $N$=14, $L$=10, $cdist$=1.5, $kmax$=200, three years monitoring.

## Acknowledgements

I would like to thank Dr. Tomasz Kaczynski, professor at Sherbrooke university and Dr. Odile Marcotte, Deputy Director at the Centre de Recherches Mathématiques, for their comments and suggestions during the writing of this article. I should also acknowledge the support provided by my administration which gave me all its trust and the budgets needed to develop and promote this new statistical vision. Mention must also be made of my colleagues who are responsible for implementing the algorithms in industrial systems, in particular engineer-researcher Michel Gauvin, who was the first to implement these algorithms in a monitoring system which is currently operational on our equipment.